# Le théorème des accroissements finis comme question curriculaire

Jean-Pierre Bourgade

Institut National Polytechnique de Toulouse, France

**Abstract**. Beyond the difficulty to give true practical motives to introduce a theorem, the didactic transposition of scholarly knowledge in a school setting often leads to freezing a technical tool into a theorem given under minimal assumptions. Probably whole categories of exercises have no other justification than establishing the necessity of introducing a theorem *under minimal assumptions* in a curriculum. The mean value theorem represents a paradigmatic situation, showing not only the expectations of the curriculum writers, but also part of the school unconscious that manifests itself in the very notion of *rigor*.

**Resumen**. Más allá de la dificultad para proponer verdaderas motivaciones prácticas para la introducción de un teorema, la transposición didáctica de un saber sabio en el ámbito escolar lleva a menudo a fijar una herramienta técnica en un teorema enunciado bajo hipótesis mínimas. Ciertas categorías de ejercicios probablemente no tienen más justificación que la de establecer la necesidad de introducir en el currículo un teorema *bajo sus hipótesis mínimas*. El teorema del valor medio constituye una situación paradigmática que permite observar no solo las expectativas de los autores del currículo sino también una parte del inconsciente escolar que se manifiesta en la noción de *rigor*.

**Résumé.** Outre la difficulté à proposer de réelles motivations pratiques pour l'introduction d'un théorème, la transposition didactique d'un savoir savant dans le cadre scolaire conduit souvent à figer un outil technique en théorème énoncé sous des hypothèses minimales. Des catégories entières d'exercices n'ont probablement pas d'autre raison d'être que de justifier la présence au programme d'un théorème *sous ses hypothèses minimales*. Le théorème des accroissements finis constitue une situation paradigmatique pour observer non seulement les attentes des rédacteurs d'un programme d'enseignement, mais également une partie de l'inconscient d'école qui se manifeste dans la notion même de *rigueur*.

Jean-Pierre Bourgade

## 1. La rigueur et l'invention

Les mathématiques sont souvent perçues comme un domaine aride où règne « la rigueur », c'est-à-dire une police de la pensée et de l'expression qui rejette les arguments invalides, les expressions floues, les calculs inexacts. Cette perception n'est pas nécessairement celle qu'en ont les mathématiciens eux-mêmes qui peuvent mettre en avant le rôle de l'imagination, de l'inventivité, du flou, du hasard dans la construction de leurs idées. La notion de théorème est, elle aussi, fortement attachée à l'idée qu'on se fait des mathématiques et les premiers théorèmes rencontrés structurent fortement l'image que l'on se fait de cette notion : le théorème de Pythagore *s'applique* à des triangles pour autant que certaines *hypothèses* soient satisfaites. En général, le théorème de Pythagore est énoncé sous ses *hypothèses minimales* (on considère un triangle rectangle, mais on ne précise pas s'il est isocèle ou non – pourtant le théorème reste vrai sous l'hypothèse plus forte d'un triangle rectangle isocèle) et tout le monde en est satisfait car cela permet de l'appliquer à un plus grand nombre de situations que si les hypothèses étaient plus difficiles à satisfaire (il est plus rare de rencontrer un triangle rectangle isocèle qu'un triangle rectangle quelconque). Une autre caractéristique des théorèmes est qu'ils sont démontrés *rigoureusement*, et que leur utilisation (correcte) permet à son tour de parvenir à des *conclusions* certaines de façon rigoureuse. Pourtant, un théorème n'est qu'un point saillant dans une branche de l'activité humaine et la réalité conduit à des adaptations qui ne se fondent pas aisément dans le moule de cette rigueur schématique. La notion de triangle rectangle n'est qu'un modèle, on peut espérer appliquer le théorème de Pythagore à des triangles à peu près rectangles et estimer l'erreur commise. Ce type de pratiques est pratiquement inexistant à l'école, y compris dans les lieux où se forment les futurs ingénieurs français (principalement les classes préparatoires aux grandes écoles, CPGE, qui accueillent des étudiants âgés de 18 à 20 ans environ avant leur accession sur concours aux écoles d'ingénieur), où la rigueur reste perçue comme un canon, contraignant pour certains, plaisant pour d'autres, mais en tous cas jamais remis en question.

À l'instar du théorème de Pythagore, le théorème des accroissements finis (TAF) est un outil central, à la fois théoriquement et technologiquement puisqu'il permet non seulement de démontrer d'autres théorèmes mais aussi de fonder des techniques de calcul approché en contrôlant leur précision. Ce type de résultat offre les conditions idéales pour observer *in vivo* la lutte qui oppose les idéaltypes de la rigueur et de l'imagination. Le TAF est actuellement au programme des CPGE scientifiques françaises, présenté sous ses hypothèses minimales. Cette situation présente cependant une délicatesse qui n'apparaît pas dans le cas du théorème de Pythagore : il est facile de trouver des triangles rectangles non isocèles, c'est-à-dire des triangles qui satisfont les hypothèses minimales du théorème de Pythagore *sans plus*, alors qu'il est pratiquement impossible de rencontrer (dans des contextes pertinents d'un point de vue pratique) des objets (en l'occurrence des fonctions) qui remplissent *sans plus* les conditions de validité minimales du TAF. Cette situation, assez courante en analyse, conduit à une tension forte entre deux contraintes liées à l'enseignement en CPGE : le désir de rigueur conduit à énoncer le théorème sous ses hypothèses minimales et à le démontrer dans ce cadre, de façon à maximiser la généralité d'application du théorème ; par ailleurs, les applications réalistes du théorème – ses raisons d'être, au sens de la théorie anthropologique du didactique (TAD) – ne font appel qu'à des conditions de validité moins larges et ne légitiment donc pas pleinement la nécessité d'énoncer le TAF sous ses hypothèses minimales. Cela conduit à deux types de pratiques opposées, qui ne font qu'exprimer de façons différentes la même tension. D'une part apparaissent des exercices « théoriques » dont l'objectif est de faire démontrer un résultat général, valable pour une classe d'objets (de fonctions) très vaste et caractérisée de façon *ad hoc* par le fait que ses objets satisfont *sans plus* les conditions minimales de validité du TAF. Ils ne font que déplacer le problème de la recherche des raisons d'être des hypothèses minimales du TAF : où sont les raisons d'être de ces exercices ? D'autre part, les exercices qui correspondent aux situations réalistes d'utilisation du TAF peuvent être traités sans faire appel aux conditions minimales de validité de ce théorème, mais la présence explicite au programme des CPGE du théorème sous ses hypothèses minimales contraint les enseignants, les étudiants, les rédacteurs de manuels à calquer leurs attentes sur le modèle suivant : « qui peut le plus peut le moins, mais il faudra bien veiller à s'assurer dans chaque cas que qui peut le plus peut effectivement le moins ». La situation est comparable à celle d'un monde où la plupart des triangles seraient rectangles isocèles, et où on persisterait à énoncer le théorème de Pythagore pour les triangles rectangles quelconques tout en demandant aux étudiants de bien vérifier en chaque cas que tel triangle rectangle isocèle est bien, en outre, un triangle rectangle (non nécessairement isocèle). Dans le cas du TAF, cette exigence est souvent vécue comme une exigence de rigueur, par les enseignants comme par les étudiants.

Une autre utilité est avancée pour justifier la présence au programme du TAF sous ses hypothèses minimales : il permet de démontrer rigoureusement d'autres théorèmes fondamentaux. De nombreuses voix se sont élevées pour critiquer ce point de vue et plusieurs propositions ont été faites pour donner d'autres démonstrations plus directes des résultats en question, sans faire appel au TAF. On arrive donc à une situation paradoxale, où quasiment rien ne plaide pour la persistance au programme d'un théorème énoncé sous ses hypothèses minimales, mais où il demeure malgré tout, défendu avec force par beaucoup d'enseignants et de noosphériens. Cette survie d'un savoir moribond montre bien qu'il n'y a pas, en matière curriculaire, de force intrinsèque de l'idée vraie, ou bien qu'à cette force s'opposent d'autres forces de réaction (dans tous les sens du terme), relevant d'un registre plus sociologique que logique.

L'objectif de cette communication est d'essayer de rendre compte de ces forces, afin de produire une explication de la persistance dans le curriculum d'un théorème faiblement motivé. L'étude historique de l'émergence du TAF vu comme un résultat théorique qui permettait de fonder l'analyse, non seulement en tant que branche des mathématiques, mais surtout en tant que champ scientifique fortement intégré au sens de la théorie des champs de Pierre Bourdieu (2001), l'observation de la transposition progressive de ce savoir savant en un savoir scolaire, permettent de comprendre la place symbolique qui est faite à la notion de rigueur dans un modèle d'enseignement centré sur le monumentalisme.

En particulier, nous conjecturons que le fait que le TAF donne lieu à des séries d'exercices permettant un traitement algorithmique du raisonnement fournit à l'enseignant des moyens techniques d'évaluer « la rigueur » de ses étudiants et, en retour, *d'imposer comme légitime* une certaine vision de la rigueur.

## 2. Lente émergence d'un résultat « théorique »

Au XVII$^e$ siècle, l'un des opposants les plus ardents à l'analyse naissante, Michel Rolle, élabore une technique de localisation et comptage des racines réelles de polynômes, la « méthode des cascades ». Cette méthode affirme que les racines d'un polynôme sont séparées par les racines de son polynôme dérivé (sa « cascade ») et l'itération de ce procédé permet de localiser grossièrement les racines. La praxéologie construite par M. Rolle dans son *Traité d'Algèbre* (1690) est vivement critiquée faute d'une technologie complète, ce qui conduit M. Rolle à publier une *Démonstration d'une méthode pour résoudre les égalitez de tous les degrez* (1691), démonstration purement algébrique d'un résultat qui ne porte que sur les polynômes (pour une exposition plus détaillée de la méthode des cascades voir David E. Smith, 1984, pp. 253-260 et Julius Shain, 1937).

Il faut attendre le XIX$^e$ siècle pour voir renaître le résultat de M. Rolle sous la forme d'un résultat d'analyse. L'émergence du théorème de Rolle sous ses hypothèses minimales[1] apparaît comme une réponse à un problème différent du type de tâches considéré par M. Rolle : comment justifier aussi complètement que possible (sachant que ce possible n'était pas connu avant d'être réalisé) les pratiques des mathématiciens en analyse au moyen de raisonnements fondés sur des résultats eux-mêmes déjà démontrés ou admis pour vrais (axiomes) ? Cet aspect de l'histoire de l'analyse montre le lien consubstantiel qui se crée progressivement entre l'exigence de rigueur et la minimalité des hypothèses d'un théorème. On est pourtant encore bien loin de la rigueur pour la rigueur

---

1. Si $f$ est une fonction continue sur [$a$,$b$] et dérivable sur ]$a$, $b$[, telle que $f(a) = f(b)$, alors il existe un point $c$ situé dans l'intervalle ]$a$,$b$[ en lequel $f'(c) = 0$. En particulier, si une fonction possède deux racines, elles sont séparées par une racine de sa dérivée, ce qui est au cœur de la méthode des cascades de Rolle.

comme en témoigne Jules Tannery cité par Hélène Gispert-Sambaz (1982) :

> On peut raisonner fort bien et fort longtemps sans avancer d'un pas, la rigueur n'empêche pas un raisonnement d'être inutile […]. S'imagine-t-on, par exemple les inventeurs du calcul différentiel s'acharnant avant d'aller plus loin, sur les notions de dérivée et d'intégrale définie ? Ne valait-il pas mieux montrer la fécondité de ces notions, dont l'importance justifie le soin qu'on a mis à les éclaircir ? Cette révision même, qu'on a faite de notre temps, l'aurait-on entreprise, sans les questions que l'étude des fonctions et particulièrement des séries trigonométriques a posées d'une manière inévitable ? (p. 36)

On voit que la rigueur n'a émergé en analyse que poussée, d'une manière inévitable, par l'étude de questions. Dans un deuxième temps, sans doute en lien avec l'exaspération des élèves allemands de Karl Weierstrass de n'être pas compris par leurs collègues (Gispert-Sambaz, 1982, p. 33), puis sous l'impulsion de Ulisse Dini, qui voit dans l'approche de K. Weierstrass une confirmation de ses propres doutes (Gispert-Sambaz, 1982, p. 30), et aussi de Gaston Darboux (Gispert-Sambaz, 1982, p. 67 sq.) qui constate l'indigence de la plupart des traités d'analyse français, paraissent en Europe des traités d'analyse reprenant les idées de K. Weierstrass (dus à U. Dini, Camille Jordan, Carl Harnack, Paul du Bois Reymond, Otto Stoltz, etc.). À compter de ce moment, la raison d'être de ce surcroît de « rigueur » n'est plus seulement de donner des réponses à des questions pratiques (sur les séries de Fourier, etc.), mais aussi et surtout de proposer un exposé clair et *self-consistent* de l'analyse à l'usage des nouveaux entrants dans le champ des mathématiques, ou dans le champ récemment renouvelé de l'analyse, pour lequel les nouveaux entrants peuvent être des mathématiciens déjà aguerris [2]. Si la volonté de faire école n'est pas toujours explicite, et si elle s'exprime surtout comme un désir de ménager une inter-intelligibilité pour les analystes, elle apparaît en creux chez les élèves allemands de K. Weierstrass et de façon plus explicite en Italie où la publication du traité de U. Dini fait de ce dernier le chef de file de la nouvelle école d'analyse italienne (Gispert-Sambaz, 1982, pp. 31-32). Le besoin de *faire corps* conduit à l'émergence de formes de justification communes qui autorisent la compréhension et donc aussi l'évaluation mutuelles des mathématiciens – même si cette tendance n'est pas homogène en Europe, elle est plus vive en Allemagne et en Italie qu'en France par exemple, où G. Darboux reste isolé comme le souligne H. Gispert-Sambaz (1982, pp. 67 *sq*.). Le théorème de Rolle et le TAF sont au cœur de ce processus puisque la *création* progressive du théorème de Rolle *sous ses hypothèses minimales* commence avec la démonstration du TAF au moyen du théorème de Rolle par Pierre-Ossian Bonnet (1844 ; voir aussi H. Gispert-Sambaz, 1982, p. 23) puis fait l'objet d'une intense discussion qui vise à effacer toute incertitude sur la démarche (imparfaite) de P.-O. Bonnet : G. Darboux, puis U. Dini améliorent la preuve de P.-O. Bonnet et le mouvement s'achève lorsque Hermann Schwarz déduit du TAF le principe de Lagrange [3] pour la première fois sans supposer la continuité de la dérivée. C'est donc par un hasard historique que le théorème de Rolle et le TAF se trouvent placés au cœur de « *la seule* exposition rigoureuse et simple du calcul différentiel » (Darboux, cité par Gispert-Sambaz, 1982, p. 24 ; c'est nous soulignons). Cette illusion (puisque c'en est une) perdure longtemps et explique largement la persistance de ces théorèmes aux programmes de l'enseignement secondaire et supérieur.

---

2. Sur ces questions, voir P. Bourdieu (2001).

Les programmes de 1905 de terminale C (élèves de 17 ans, dominante scientifique) précisent sans ambiguïté que « le professeur laissera de côté toutes les questions subtiles que soulève une exposition *rigoureuse* de la théorie des dérivées » ainsi que le note Jean-Pierre Daubelcour (2009 ; c'est nous qui soulignons). Le programme de 1962 de cette même classe, qui cherche plus de rigueur dans l'exposé de la théorie de la dérivation, inscrit pour la première fois au programme le théorème de Rolle (admis) et le TAF (démontré à partir du théorème de Rolle et dont on déduit le principe de Lagrange). L'arrivée du TAF coïncide avec la création d'une rubrique « Analyse » dans le programme de 1962 : auparavant, l'étude des fonctions dérivables était intégrée dans la partie « Algèbre » du programme ; on peut voir ici une lointaine réplique du tremblement de

---

3. Le principe de Lagrange affirme que toute fonction dont la dérivée est constamment positive sur un intervalle est croissante.

Jean-Pierre Bourgade

terre que constitue l'émergence de l'analyse comme sous-champ autonome du champ mathématique. Pourtant la raison d'être principale de l'introduction du TAF dans les ouvrages savants du XIX$^e$ siècle, le besoin d'unifier un langage et des pratiques, se sont perdus. Ici, la structure même du champ scolaire fait qu'il n'y a aucune attente du public, parents d'élèves et élèves, à cet égard : ils font une confiance toute contractuelle au bien-fondé du savoir enseigné. On assiste alors à un glissement fonctionnel qui va conduire peu à peu à vider l'exigence de rigueur de son sens initial pour conduire à la conception (implicite) actuelle de la notion de rigueur.

Ainsi, comme le souligne J.-P. Daubelcour (2009), avec l'apparition du théorème de Rolle « une phase déductive très formatrice [devient] possible […] : la démonstration du théorème des accroissements finis accompagné de son interprétation géométrique, permet la démonstration du principe de Lagrange » (p. 117). La possibilité de démontrer quelque chose devient une des raisons d'être de l'inscription du TAF au programme.

Par la suite le TAF disparaît du programme avec la réforme de 1971, dite des « mathématiques modernes », le principe de Lagrange étant admis dès la classe de 1$^{re}$. Pourtant, le TAF réapparaît en 1983, mais afin d'en déduire l'inégalité des accroissements finis (IAF) et non le principe de Lagrange qui reste admis en classe de 1$^{re}$. Enfin, à partir de 1986, le TAF disparaît définitivement, même si l'IAF perdure au programme jusqu'à la fin du XX$^e$ siècle. On voit que l'objectif est de recentrer le programme sur les aspects quantitatifs : « majorations, encadrements, vitesse de convergence, approximation à une précision donnée » (MEN[4], 1986), quitte à minorer la part de justification théorique des techniques en question. Le programme de 1983 précisait déjà que « dans les énoncés et les démonstrations on continuera de se placer dans des hypothèses de bonne sécurité sans en rechercher de plus fines ». Le but n'est donc plus de « fonder l'analyse » et les exercices proposés dans les manuels correspondent à ces objectifs. Le désir fondationnel deviendra l'apanage de l'enseignement supérieur : le théorème de Rolle, le TAF ainsi que leur démonstration sont au programme des CPGE scientifiques en 1984 et leur utilisation dans la plupart des manuels reproduit l'ambiguïté essentielle du rôle que joue le TAF. On retrouve, d'une part, la tentation théoriciste avec des exercices faisant appel aux hypothèses minimales ; et, d'autre part, des exercices mettant en valeur les « aspects quantitatifs », du type « majorer, minorer, estimer une erreur ».

---

4. Le sigle MEN fait référence au ministère de l'Éducation nationale.

## 3. Le TAF et ses usages en CPGE

Le TAF est généralement énoncé pour des fonctions définies sur un intervalle borné de la forme [$a$, $b$]. On suppose qu'une certaine fonction $f$ est continue sur [$a$, $b$] et dérivable sur ]$a$, $b$[. Le TAF affirme qu'il existe alors un point noté $c$, contenu dans l'intervalle ]$a$, $b$[, tel que $f'(c) = (f(a) - f(b))/(b - a)$. Le programme indique l'interprétation géométrique suivante : la pente moyenne d'une courbe entre ses extrémités est égale à la pente de la courbe en au moins un point ; ainsi que l'interprétation cinématique suivante : la vitesse moyenne d'un mobile entre deux instants est égale à la vitesse instantanée du mobile à un instant au moins. Les deux interprétations sont trompeuses car elles ne correspondent pas réellement aux hypothèses minimales du théorème : ainsi, les variations de vitesses d'un mobile ne peuvent être trop fortes sous peine de destruction du mobile. Le théorème sous ses hypothèses minimales n'est donc pas adapté aux tâches de modélisation. Pourtant, le TAF fait l'objet du commentaire suivant dans les programmes des CPGE Mathématiques, Physique et Sciences de l'Ingénieur (MPSI) publiés en 2003 : « Pour le théorème de Rolle, l'égalité et l'inégalité des accroissements finis, ainsi que pour la caractérisation des fonctions monotones, on suppose $f$ continue sur [$a$, $b$] et dérivable sur ]$a$, $b$[. » (MEN, 2003)

L'émergence historique du TAF sous ses hypothèses minimales permet de comprendre l'insistance du programme sur ce point : le but est de permettre « enfin » la démonstration du principe de Lagrange, admis et utilisé depuis la classe de 1$^{re}$. On fait généralement deux types d'objections à un tel point de vue. La première, classique et avancée par exemple par Jean Dieudonné (1969, p. 148), Ralph Boas (1981), Thomas W. Tucker (1997), Jean-Louis Ovaert et Jean-Luc Verley (1983, p. 147), Christian Houzel (1996, p. 8), consiste à remarquer que, de toutes

Jean-Pierre Bourgade

façons, l'*égalité* des accroissements finis n'apporte qu'une précision illusoire puisque la démonstration du TAF n'est pas constructive et ne permet pas de localiser le point $c$. L'IAF suffit à la démonstration du principe de Lagrange. Certains, comme T. W. Tucker (1997), prolongent la critique en remarquant qu'il est choquant de « assum[e] the nonobvious [c.-à-d. l'IAF] to prove the obvious [le principe de Lagrange] » (p. 231), critique que Henri Lombardi (1999) reprend avant de souligner que ce type d'arguments passe à côté de l'essentiel : « on ne devrait jamais démontrer le théorème des accroissements finis […] en tant que preuve d'une chose évidente, mais en tant que vérification de l'adéquation d'un modèle mathématique à la réalité qu'il veut représenter » (p. 56). Bien sûr, ce point de vue qui pourrait donner une légitimité au TAF n'est pas adopté en CPGE, où la question de la modélisation n'est jamais abordée frontalement malgré les vœux pieux que l'on trouve ici ou là dans les programmes et qui sont aussitôt contredits par l'utilisation très significative des notions d'*interprétation* et d'*application*. En outre, H. Lombardi souligne à juste titre que, du point de vue de la modélisation, l'hypothèse de dérivabilité est beaucoup trop générale (essentiellement : la vitesse d'un mobile ne peut pas être une dérivée « générique ») et qu'on pourrait donc se contenter de versions beaucoup plus faibles du TAF (Lombardi, 1999, p. 63).

## 4. Illustration de l'applicationnisme

La justification de la présence du TAF par la nécessité de démontrer le principe de Lagrange se réduit au fond au propos de J.-P. Daubelcour (2009) sur les « phases déductives très formatrices » : si le programme impose l'étude du TAF sous ses hypothèses minimales, c'est probablement au nom de la mise en œuvre de raisonnements « formateurs » que cette situation rend possibles. Le programme place la classe dans une problématique interventionniste : on est conduit à explorer l'ensemble $\{\Pi \ / \ \Im \ (\wp, \Pi, U)\}$ des projets $\Pi$ (c'est-à-dire, finalement, des études d'exercices) de l'instance $U$ (les élèves) pour lesquels il est utile ou indispensable de disposer d'une praxéologie $\wp$ qui reste à préciser mais qui semble condamnée par le programme à tourner autour du TAF sous ses hypothèses minimales, ce qui pose des difficultés importantes.

D'une part, il faut trouver des exercices où le théorème sous ses hypothèses minimales est utile, voire indispensable. On assiste alors à une prolifération d'exercices où l'utilité du théorème n'est même plus à prouver puisque les situations sont créées pour éviter le problème de la recherche des raisons d'être. Ce sont tous les exercices dont l'énoncé commence par « Soit une fonction continue sur $[a, b]$, dérivable sur $]a, b[$ telle que... ». On trouvera une série d'exemples dans le manuel de Jean-Marie Monier (1999) qui explicite ainsi les « raisons d'être » du théorème de Rolle dans le programme des classes préparatoires :

> Pour établir une propriété du type « il existe $c \in ]a, b[$ tel que... », la fin de la propriété faisant intervenir une dérivée, on peut essayer :
> – d'appliquer le théorème de Rolle ou le théorème des accroissements finis, une ou plusieurs fois […],
> – de construire un réel A et une fonction auxiliaire $\varphi$ en s'inspirant de la preuve du théorème des accroissements finis […]. (p. 141)

D'autre part, les exercices où le théorème pourrait s'appliquer sans qu'il soit nécessaire de l'invoquer sous ses hypothèses minimales sont peu ou prou exclus *a priori* du fait qu'ils ne font pas apparaître comme particulièrement utile ni indispensable l'utilisation du théorème *sous ses hypothèses minimales*. Or les situations concrètes qui donnent ses raisons d'être à un résultat comme l'IAF (étude de systèmes dynamiques par exemple) tombent dans cette catégorie. Bien entendu, en pratique, on contourne la difficulté en étudiant ce type d'exercices malgré tout, mais l'injonction du commentaire du programme de MPSI conduit à faire comme si le théorème, sous ses hypothèses minimales, était utile ou indispensable dans des cas où il ne l'est peut-être pas.

Il apparaît donc qu'un théorème comme le théorème de Rolle, autrefois fonctionnellement motivé, a perdu ses raisons d'être pratiques pour devenir un lemme qui trouve son insertion naturelle dans une construction axiomatique possible de l'analyse moderne, puis un point de départ pour de possibles « phases déductives formatrices » pour les étudiants. Il reste à préciser la fonction de ce rôle que jouent actuellement

le TAF et le théorème de Rolle en CPGE. L'analyse menée par Yves Chevallard (1991) dans un contexte différent peut donner une piste :

> [I]l existe un contrat-type à propos du « comment travailler ? ». L'élève est en effet censé « apprendre ses leçons et faire ses exercices » […] le professeur sera, lui, en principe inattaquable dès lors qu'il aura mis entre les mains de l'élève ce qui permettra à celui-ci d'accomplir son devoir d'élève : des leçons à apprendre, des exercices à faire. (pp. 176-177)

Si on entre dans les détails :

> [I]l existe des cas où apprenabilité et faisabilité sont *ultra-légitimes* : il s'agit des leçons pour lesquelles « apprendre » signifie « apprendre un texte », par exemple « apprendre des définitions » ; et des exercices pour lesquels « faire » signifie « faire usage d'un algorithme (de calcul) ». (Chevallard, 1991, p. 177)

À l'instar de la notion de distance en classe de $4^e$ dans la réforme des mathématiques modernes – élèves de 13-14 ans – mais aussi des développements limités, ou des définitions axiomatisées des normes, des produits scalaires, ou encore des définitions « epsiloniques » des limites dans l'enseignement supérieur, les théorèmes comme le théorème de Rolle et le TAF constituent des morceaux de choix. En effet, ils proposent un texte à apprendre précisément (l'énoncé du théorème et en particulier ses hypothèses, forcément minimales puisqu'un accord intersubjectif doit se faire sur ce qui est à apprendre), à savoir restituer, par exemple dans le cadre d'une « question de cours ». Ils donnent également l'occasion d'une pratique algorithmique qu'illustrent bien les conseils de J.-M. Monier cités précédemment : dans un sens étendu du mot, il y a calcul dès lors qu'on doit appliquer le théorème de Rolle ou le TAF puisqu'une certaine procédure réglée, assortie d'étapes systématiques (vérifier que la fonction est bien continue sur un segment, dérivable à l'intérieur du segment, vérifier l'égalité des valeurs de la fonction aux bornes du segment – pour le théorème de Rolle –, calculer le taux d'accroissement de la fonction aux bornes de l'intervalle – pour le TAF) est nécessaire pour appliquer le théorème, si bien que même si le résultat visé n'est pas atteint, l'expert pourra reconnaître le respect scrupuleux d'un algorithme chez un élève s'étant conformé au contrat, ou au contraire s'indigner que tel élève n'ait « même pas su faire ça » alors qu'il « suffit d'apprendre par cœur ».

Les conséquences les plus évidentes de ce phénomène sont lisibles dans les copies d'étudiants ou sur les forums (voir figure 1).

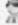

*Figure 1*. Exemple d'utilisation du TAF [5].

La sacralisation des conditions de validité du théorème est patente : il est indispensable de préciser que la fonction exponentielle est dérivable sur l'ouvert ]0, 1[. Quelle fonction remplit une telle insistance dans la rédaction d'une solution (dont on peut supposer qu'elle sera reproduite dans un devoir à rendre, en examen, etc.) ? L'objectif inconscient est probablement d'envoyer un signal au lecteur du forum tout comme, le cas échéant, à l'enseignant-évaluateur : « Je connais les hypothèses du théorème et vous ne pouvez donc pas me sanctionner sur ce point ». Voici donc une première conséquence : les étudiants sont conduits à écrire, à souligner, à mettre en avant des propos qui n'ont pas de pertinence en soi. Par ailleurs, soulignons sur cet exemple l'algorithmisation en œuvre dans la présentation du raisonnement : l'ordre et la nature des étapes sont parfaitement codifiés.

---

5. Exemple accessible à l'adresse http://www.ilemaths.net/forum-sujet-478832.html.

Jean-Pierre Bourgade

Si l'on se demande d'où peut venir une telle propension des étudiants, on trouvera la réponse dans l'exemple de corrections d'exercices donné en figure 2.

Soit $f(y) = \ln(y)$. La fonction $f$ est dérivable sur $]0, +\infty[$ et $f'(y) = \frac{1}{y}$, qui est décroissante sur $\mathbb{R}^{+*}$. Soit $x > 0$. On applique le théorème des accroissements finis sur $[x, x+1]$ ($f$ y est continue, et dérivable sur $]x, x+1[$). Il existe $c \in ]x, x+1[$ tel que $f(x+1) - f(x) = (x+1-x)f'(c) = f'(c)$. Donc
$$\frac{1}{x+1} < \ln(x+1) - \ln(x) < \frac{1}{x}.$$

Figure 2. Un exemple de correction d'exercice par un enseignant [6].

Dans l'esprit même de l'enseignant – et donc probablement dans l'impensé professoral en général puisque ce qui est visé ici n'est pas la pratique de tel ou tel professeur, mais plutôt ce qui s'exprime à travers elle d'un *habitus* professionnel –, l'important, ce sont « les hypothèses », qu'il s'agit de bien « vérifier » avant « d'appliquer » le théorème. Pourquoi un enseignant, connaisseur et praticien des mathématiques, se sent-il obligé de préciser que la fonction logarithme est dérivable sur l'ouvert $]x, x+1[$ ? Pourquoi ne se contente-t-il pas de dire qu'elle est indéfiniment dérivable sur $]0, +\infty[$ et qu'elle vérifie donc la propriété des accroissements finis ? Sans doute parce qu'il n'existe rien de tel qu'une propriété des accroissements finis dans la culture scolaire. Ce qui existe, et fortement, c'est un théorème, assorti de ses « hypothèses ». Nous avons vu comment le théorème s'est peu à peu figé avec ses hypothèses minimales : à ce stade, il devient indispensable d'invoquer ces hypothèses au moment d'utiliser le théorème.

Il est difficile de ne pas voir ici une illustration de « l'enfermement dans les questions évaluables » qu'ont mis en évidence Verónica E. Parra et María R. Otero (2011) dans le cadre des universités argentines. On peut en effet se demander si l'une des raisons d'être principales de la présence du TAF et des exercices que nous venons d'évoquer n'est pas qu'ils constituent une réponse à la question suivante : « comment évaluer la capacité des étudiants à appliquer un théorème donné ? ». Comme le notent Berta Barquero, Marianna Bosch et Josep Gascón (2007), « la necesidad de evaluar la eficacia de los procesos de enseñanza y aprendizaje de las matemáticas en las instituciones didácticas […] tiende a provocar una diferenciación y autonomización interna del corpus enseñado, así como una mayor algoritmización del mismo con la consiguiente pérdida de sentido funcional del saber enseñado ». Ainsi, l'IAF disparaît au bénéfice du TAF qui se prête mieux à une algorithmisation et autorise la fabrication d'un type d'exercices *ad hoc* qui permet d'évaluer en retour la maîtrise de l'algorithme.

---

6. http://armana.perso.math.cnrs.fr/enseignement-archives/200708-lm110/lm110-td4-correction.pdf.

L'esprit applicationniste joue d'ailleurs un rôle important de ce point de vue puisqu'il représente une réalisation du processus d'algorithmisation de la pratique mathématique elle-même. En effet, comme tous les formalismes, il conduit à ignorer les raisons d'être des théorèmes introduits dans un curriculum en créant de façon artificielle des « situations de monopole technico-technologique » (là où plusieurs techniques pourraient être efficaces, il n'en maintient qu'une seule ; là où plusieurs justifications des techniques seraient envisageables, une seule survit) afin de promouvoir une vision algorithmique de la pratique mathématique : à toute situation problématique correspond une technique qu'il convient (et qu'on postule qu'il est possible) d'identifier par simple exploration systématique du bagage technique acquis en cours, puis qu'il convient d'appliquer en respectant un protocole fixé une fois pour toutes (vérifier que les conditions de validité du théorème sont respectées, etc.) avant de tirer la conclusion attendue et souvent explicitée dans l'énoncé même du problème. On notera que cette algorithmisation contribue à favoriser une approche rétrocognitive de l'enseignement.

Même si nous ne l'utilisons pas dans le contexte de la modélisation mathématique (pour les sciences expérimentales), le concept d'applicationnisme tel que nous le détaillons ici n'est probablement qu'un cas particulier de celui qu'ont proposé B. Barquero, M. Bosch et J. Gascón (2010) – un cas d'applicationnisme interne par opposition à l'applicationnisme externe qui mettrait en jeu une autre discipline à laquelle on « appliquerait » les mathématiques. Les caractéristiques de l'applicationnisme sont vérifiées ici aussi : « Las matemáticas se mantienen independientes de las otras disciplinas » (ce qui ferait office ici « d'autre discipline » c'est, par exemple dans le cas du théorème de Rolle, la question de la localisation des racines d'un polynôme) ; « La enseñanza de las matemáticas sigue la lógica deductivista » ;

« proliferación de cuestiones aisladas » ; et surtout « La enseñanza de las herramientas matemáticas básicas siempre es anterior a l'estudio de su aplicación ».

## 5. Évaluer le sacré

Paradoxalement, c'est lorsqu'il s'éloigne le plus d'une certaine vérité épistémologique concernant les praxéologies mathématiques effectives que l'« applicationnisme » contribue à renforcer une image largement partagée de ce que sont censés être le savoir et la pratique mathématiques. Par un évitement systématique de toutes les situations de modélisation (intra- ou extra-mathématique), l'enseignement classique des mathématiques organise une amnésie collective des laïcs comme des clercs (de la société, des enseignants, des étudiants, et même des chercheurs en mathématiques qui, selon une des ruses les plus fréquentes de l'intelligence, véhiculent des représentations de leur activité qui peuvent être largement en conflit avec la réalité de leur pratique) à propos de ce que faire des mathématiques veut dire : il ne s'agit plus de chercher à répondre « *by any means necessary* » à des questions qui se posent, et donc notamment, mais entre autres, par des moyens mathématiques, mais « d'apprendre les mathématiques », souvent au moyen d'un entraînement épuisant et privé de sens.

Ce qui naît alors, c'est une illusion sur ce qui fait le cœur de la pratique mathématique : la « rigueur mathématique » (Chevallard, 2011). À l'opposé des incertitudes liées à la pratique de la modélisation, dans laquelle on n'est jamais sûr d'avoir trouvé le bon modèle, ou plutôt dans laquelle on est assuré de n'avoir construit qu'un modèle, la pratique enseignante fondée sur la répétition d'exercices semblables articulés sur l'application de théorèmes, de techniques, etc., stéréotypés, donne corps au mythe d'une rigueur mathématique qui ferait consensus et qui se résumerait à :

– l'univocité des techniques mises en œuvre et l'unicité du résultat obtenu par leur application,
– la nécessité, au-delà de l'emploi de certaines techniques, d'organiser l'utilisation de ces techniques dans un discours lui-même fortement codifié, la codification du discours étant un des éléments indispensables pour garantir la rigueur du « raisonnement »,

– la spécificité du raisonnement mathématique (distingué du raisonnement du physicien, du chimiste, du sociologue, etc.), dépourvu d'ambiguïté, susceptible de vérifications quasi-automatiques,

– la transparence, la naturalité, l'intemporalité et l'ubiquité des principes de codification du raisonnement mathématique rigoureux.

On peut opposer à ces principes :

– qu'une même tâche mathématique (y compris un calcul) peut être réalisée de multiples façons (par exemple, un étudiant de terminale pourrait établir l'inégalité $e^x \geq x+1$ en étudiant la fonction $f(x) = e^x - x - 1$, et donc sans faire appel au TAF),

– que, contre une certaine fétichisation du langage mathématique, il est possible de proposer des démonstrations valables en langue vulgaire, voire au moyen de dessins, comme le fait Roger B. Nelsen (1993),

– que tout raisonnement, toute pensée, toute expression est par essence potentiellement ambiguë et que, comme le souligne Ludwig Wittgenstein (2005), toute règle est en réalité interprétée et on ne saurait, sans risque de régression à l'infini, proposer de règle pour l'usage de la règle,

– que les canons mêmes de la rigueur telle que nous la connaissons sont très variables selon les époques (y compris chez Augustin Cauchy), les lieux, et même les moments (voir le relâchement de l'enseignant redevenu mathématicien qui s'autorise évidemment les raccourcis mêmes qu'il interdit à ses étudiants).

Néanmoins, cette vision de la rigueur mathématique et sa sacralisation sont sans doute à compter au nombre des contraintes majeures qui pèsent sur la possibilité d'un changement de paradigme scolaire. On peut y voir une des raisons d'être principales de la survie du théorème de Rolle dans le curriculum, auprès des définitions de la convergence (des suites, etc.), des définitions axiomatiques diverses (norme, produit scalaire, etc.), des théorèmes de convergence (des suites, des séries, des intégrales généralisées), etc. : il propose une situation idéaltypique où peut s'exprimer la virtuosité du mathématicien rigoureux : connaissance (par cœur) de définitions (ou d'hypothèses, etc.), vérification d'hypothèses (ou

Jean-Pierre Bourgade

d'axiomes), application (à une situation particulière), algorithmisation possible de l'application (dans un cadre relativement large d'exercices possibles), etc. L'algorithmisation va de pair avec la possible codification de la rédaction (avec des étapes imposées, etc.) qui autorise en retour une évaluation facilitée du sacré : la rigueur qui est censée être l'essence de la pratique mathématique est à la fois l'objet d'une sacralisation (Chevallard, 2011), ce qui la rend impossible à évaluer car on ne peut évaluer que le profane, et le prototype de la compétence à acquérir – et donc à pouvoir être évaluée. Des exercices comme ceux qui abondent autour du théorème de Rolle permettent ce miracle en offrant des lieux restreints où la pratique peut être (presque) entièrement codifiée, autorisant par là même une évaluation de « la rigueur » *via* l'assimilation de la rigueur à la codification (dans les rapports des jurys des Concours Communs Polytechniques[7], l'attente de rigueur est souvent reliée explicitement d'une part, à la « vérification des hypothèses » et, d'autre part, aux « qualités de rédaction » ; on retrouve ces attentes dans la bouche des enseignants si l'on en croit par exemple cette contribution à un forum[8] sur internet : « Manquer de rigueur pour ma prof de maths c'est mal rédiger et ne pas assez développer son raisonnement »).

Du point de vue de la TAD, la construction de cours à travers des organisations didactiques *classiques*, selon la terminologie introduite par J. Gascòn (2001) et M. Bosch et J. Gascòn (2001), repose sur une valorisation exclusive des moments technologico-théorique et du travail de la technique, à l'exclusion de tout moment exploratoire notamment, et au prix d'une « trivialización de la actividad de resolución de problemas » et conduit à considérer que « la enseñanza de las matemáticas es un proceso mecánico totalmente controlable por el profesor » (Bosch & Gascòn, 2001). L'absence de moment exploratoire conduit à un *évitement* des questions liées à l'évaluation des techniques utilisées dans la réalisation des types de tâches visés et en particulier à négliger la question fondamentale de la multiplicité possible des techniques. L'exemple de la figure 1 montre bien ce refus par la mention explicite « en utilisant le théorème des accroissements finis » qui, loin de constituer une indication en vue d'aider les étudiants, représente au contraire un signe fort d'*interdiction* de toute exploration du champ des techniques possibles. Par ailleurs, la *portée* des techniques n'est pas davantage étudiée et l'insistance sur la minimalité des hypothèses relève d'un discours de type sécuritaire (« *on ne sait jamais*, il vaut mieux disposer d'un théorème de grande portée… ») qui contourne le problème en posant l'équivalence entre conditions minimales et portée maximale (ce qui peut être vrai d'un point de vue logique, mais pas d'un point de vue pratique). La monumentalisation de la rigueur conduit ainsi à faire taire les questions qui sont pourtant au cœur de la pratique mathématique, si bien que, finalement, ce n'est que dans des situations où on n'évalue au fond qu'un type de calculs dénué de raisons d'être extra-scolaires qu'on prétend appliquer la pierre de touche à la compétence par excellence, le respect de la *rigueur*.

---

7. Consultables à l'adresse suivante :
http://ccp.scei-concours.fr/sccp.php?page=cpge/rapport/rapport_accueil_cpge.html&onglet=rapports.
8. Consulter : http://forums.futura-sciences.com/mathematiques-superieur/62830-gagner-rigueur-maths.html.

## 6. Évaluer la rigueur : un rite d'institution

On est alors dans une situation paradoxale : la survalorisation symbolique de la rigueur conduit à lui ôter toute fonctionnalité. La rigueur est aimée, pratiquée, codifiée pour elle-même, en tant que « cœur de métier », alors que les mathématiques n'ont justement pas de cœur de métier mais sont intégrées de manière fonctionnelle à la pratique de multiples disciplines, notamment si l'on considère avec B. Barquero, M. Bosch et J. Gascón (2010) que « la modelización matemática debe formar parte integrante de cualquier proceso de estudio de matemáticas » puisque « la enseñanza de la modelización matemática se conviert[e] en un "sinónimo" de la enseñanza funcional de las matemáticas » (p. 554). Comme un rituel qui perdure bien que le sens des gestes se soit perdu pour les fidèles, la rigueur continue d'être « enseignée » sans que les raisons d'être de la rigueur elle-même ne soient jamais données. Les raisons habituellement fournies aux étudiants récalcitrants se résument le plus souvent à : « la rigueur en mathématiques est indispensable pour une raison toute simple : elle seule permet d'être sûr des résultats établis », comme le dit Daniel

Perrin (1997), à quoi on répondra avec T. W. Tucker (1997) que Leonhard Euler se débrouillait très bien sans la « rigueur » contemporaine, ou avec William P. Thurston (1995) que « la confiance ne vient pas de mathématiciens qui vérifient formellement les arguments formels, elle vient de mathématiciens qui étudient soigneusement et de façon critique les idées mathématiques ». L'acquisition de cette rigueur évaluable et ineffable à la fois constitue un vrai *rite d'institution* au sens que P. Bourdieu (1982) donne à ces mots : seul celui ou celle qui parvient à maîtriser le code, au prix d'une conversion d'habitus, peut être institué mathématicien (« nul n'entre ici s'il n'est rigoureux »). Se crée ainsi une frontière entre ceux et celles qui le sont (rigoureux), et les autres ; cette frontière, difficilement franchissable, prime toutes les autres différences et devient le critère par excellence de l'excellence, c'est-à-dire de la conformité à l'image que les mathématiciens se font de la pratique mathématique. En particulier, en tant que rite d'institution, l'évaluation (par examen, concours, etc.) de la « rigueur » *naturalise* un arbitraire culturel et *crée* une opposition entre ceux qui possèdent et ceux à qui manque la disposition à la rigueur (Bourdieu, 1982). Ainsi des élèves peu enclins à ajuster leur habitus aux exigences du champ scolaire, mais très inventifs, qui parviennent à trouver des réponses sans pouvoir les justifier pleinement (selon les canons en vigueur), se voient refuser l'accès à l'institution accordé à des étudiants moins riches d'idées mais plus conformes aux attentes.

Il y aurait par ailleurs tout un travail à faire sur le rôle que jouent les indications données dans les énoncés d'examens. En première analyse, on peut penser que leur présence est liée à des contraintes écologiques : besoin de la part des professeurs de présenter des sujets d'examen « intéressants », « astucieux », « originaux » (ce qui conduit souvent aux limites du hors programme et, comme nous l'avons vu, relève d'un désir de « motiver » l'enseignement plus par l'invocation d'une esthétique du sacré que par la mise en situation [9], forcément profane) et nécessité de respecter le contrat didactique en ne mettant pas les étudiants face à des types de tâches inédits un jour d'examen. L'*indication* joue probablement aussi un rôle symbolique en cela qu'elle dit sans le dire ce qu'est la distribution objective des rôles. Au professeur la maîtrise (noble) du technologico-théorique qui autorise le surplomb nécessaire à la construction d'un sujet « original », à lui les choix stratégiques (choix d'une approche, ici « en utilisant le TAF »), à l'étudiant la maîtrise (vulgaire) de la pratique de techniques qui lui permettront de traiter le sujet « original » une fois celui-ci abaissé de son piédestal au moyen d'une *indication*, véritable passe-droit qui autorise le vulgaire à accéder momentanément aux beautés aristocratiques, mais sur un mode tel qu'il n'y accède pas : à l'instar des romans populaires qui donnent à leurs lecteurs un accès provisoire et fictif aux salons de la bourgeoisie, l'indication réalise la prouesse de donner à voir tout en refusant de donner accès, renforçant par là-même les représentations que les agents (ici, les étudiants et les professeurs) se font de leurs positions relatives, mais aussi de l'activité mathématique.

---

9. Au sens de la théorie des situations didactiques.

Soulignons enfin le caractère *ineffable*, *intransmissible*, *peu objectivable* de ce qui fait pourtant l'objet de l'évaluation principale : la rigueur est indéfinissable, les attentes à son égard peuvent être très différentes d'un exercice à l'autre, d'un enseignant à l'autre, etc., mais son absence fait pourtant l'objet de sanctions sans appel. Internet regorge de forums où des étudiants échangent sur leurs difficultés en mathématiques et la lecture des contributions suivantes, qui portent sur la quête de l'inaccessible rigueur, donnent un éclairage sur ce point : nous en citerons trois ici, notés dans la suite F1, F2 et F3 [10]. On y lit déjà la confirmation du lien entre manque de rigueur et échec : « mon prof me reproche de ne pas avoir assez de rigueur, ce qui explique mon petit 12 de moyenne » (F3), « mon problème vient d'un manque de rigueur (et donc de cohérence) et c'est très pénalisant puisque même ce que je sais faire ou ce que j'ai compris ne me rapporte pas de points de ce fait... » (F2), « Lors de mes controles [*sic*] de maths je perd [*sic*] vraiment beaucoup de points car je n'ai pas assez de "rigueur" » (F1). On y trouve ensuite l'incompréhension de ce qu'est ce petit rien qui fait tout : « je n'ai pas assez de "rigueur", mais qu'est ce que c'est exactement ?? » (F1), « Je voulais donc savoir si vous sauriez (ou auriez quelques idées) comment je pourrais faire pour améliorer ma rigueur/gagner en rigueur et donc progresser ? » (F2), et enfin l'espoir que cette épineuse question puisse, comme les autres, faire l'objet d'un traitement systématique, voire automatique : « Y' as pas des trucs tout bête [*sic*] ? ou alors des automatismes a [*sic*] acquérir ? » (F3), « Des techniques, des méthodes, des astuces, des façons de travailler/réviser...Bref, comment faire quoi. » (F2). Mais surtout, de manière frappante, domine le sentiment que la réponse ne viendra pas de l'école mais s'apprend sur le tas, au petit bonheur : « j'aimerais travailler cette "rigueur" sachant que personne ne peut m'aider chez moi et que je fais deja [*sic*] plein d'exos. Vous aurez [*sic*] des conseils ?? » (F1), « la rigueur se gagne peu à peu avec l'âge et avec l'expérience dans n'importe quel domaine » (F1), « J'ai moi aussi quelques problemes [*sic*] de rigueur… Ce que j'ai trouvé de mieux, à faire, c'est [...] » (F2), « Un "truc" qui marche pas mal [...] » (F3), « J'avais le meme [*sic*] problème que toi et j'ai "trouvé" un truc tout con [...] » (F3).

---

10. F1 :   http://fr.answers.yahoo.com/question/index?qid=20081107110555AAgvoVl ;
F2 :   http://forums.futura-sciences.com/mathematiques-superieur/62830-gagner-rigueur-maths.html ; F3 : http://forums.futura-sciences.com/physique/56365-palier-a-un-manque-de-rigueur.html.

On observe là une rupture flagrante du contrat-type officiel : parce qu'il ne fournit pas les instruments permettant d'apprendre la rigueur, l'enseignant − et à travers lui, le système d'enseignement dans sa globalité − a toutes les chances de laisser les choses en l'état. La rigueur finit par gagner des titres de noblesses en apparaissant comme une qualité intrinsèque, propre à chacun, et dont certains sont dépourvus – à l'opposé des qualités qui peuvent s'acquérir par le labeur, comme « apprendre ses tables », etc. ; ce qui reconstitue en un lieu où on ne l'attendait pas, le système d'oppositions développé par P. Bourdieu et Monique de Saint Martin (1975) : noble/ignoble // pur/impur // théorique/appliqué // brillant/scolaire // sacré/profane, etc., au principe des classements scolaires. Le système d'enseignement « permet de réaliser une opération de classement scolaire tout en la masquant : il sert à la fois de relais et d'écran entre le classement d'entrée, qui est ouvertement social, et le classement de sortie, qui se veut exclusivement scolaire ». Le « vague et le flou même des qualificatifs » utilisés dans les évaluations (« peu précis », « trop vague », « pas rigoureux »), au mieux tautologiques et qui au pire montrent en creux que quelque chose manque qu'il est pourtant impossible de définir précisément, « ne véhiculant à peu près aucune information […] suffisent à témoigner que les qualités qu'ils désignent demeureraient imperceptibles et indiscernables pour qui ne posséderait déjà, à l'état pratique, les systèmes de classement qui sont inscrits dans le langage ordinaire » (Bourdieu & de Saint-Martin, 1975, p. 76). Ainsi, même dans une discipline apparemment dégagée des contingences sociales et qui offre par ailleurs un profil moins sélectif socialement que d'autres disciplines (notamment littéraires), le social n'est pas totalement absent et, refoulé dans ce qui n'est pas formalisable, le style, il est évalué par les enseignants qui, pensant rendre des verdicts objectifs sur des questions techniques, ne font peut-être, eux aussi, que sanctionner des différences d'habitus.

La question du style recoupe celle du goût : l'étudiant conforme aux attentes est d'autant plus conforme qu'il aime se conformer. Les « bons étudiants » aiment généralement les mathématiques lorsqu'ils y excellent. Chez les enseignants, le désir de « faire aimer les mathématiques », outre qu'il est contraire au principe de laïcité comme le signale Yves Chevallard (2005), s'oppose de deux façons à l'apprentissage :

– Les professeurs sont conduits à « motiver » l'introduction d'une notion à partir d'exemples « motivants » (« beaux », surprenants, spectaculaires, etc.) ; ce type de motivation interdit toute rencontre avec l'*outil* mathématique et ses raisons d'être,

– Ce type d'approches introduit au cœur de la classe un principe de vision et de division fondé sur le goût, c'est-à-dire sur des catégories sociales d'aperception du « mathématiquement beau », principe que toutes les tentatives de rationalisation de l'évaluation (barèmes détaillés jusqu'à l'absurde, etc.) ne font que refouler sur le mode de la dénégation, interdisant par là-même son objectivation et sa neutralisation relative.

Il est à noter que les étudiants ne remettent jamais en question la pertinence de l'exigence de rigueur, ni même la justesse des jugements de goût que l'enseignant ou les autres étudiants peuvent émettre à propos d'un objet mathématique : le refus de l'exigence de rigueur ne peut se faire que ponctuellement, l'allégeance étant une condition *sine qua non* de la participation au jeu – sauf à se mettre hors-jeu. L'adhésion volontaire à un système de classement dont ils sont les victimes objectives est le résultat d'une violence symbolique institutionnelle au sens que donnent à ces mots P. Bourdieu et Jean-Claude Passeron (1970) : il y a violence symbolique lorsqu'un arbitraire culturel est imposé comme légitime. On peut relire à cette lumière la lente histoire de l'émergence de la rigueur telle que nous la connaissons : d'abord *outil mathématique* et techniquement fonctionnelle, elle se transpose en un *outil de classement*, la similitude formelle (présence de symboles comme les quantificateurs, usage de tournures stylistiques récurrentes, etc.) entre les deux outils permettant la légitimation du second en déplaçant sur lui la valorisation symbolique dont le premier fait l'objet.

## Remerciements

## Références

Barquero, B., Bosch, M. & Gascón, J. (2007). La modelización matemática como instrumento de articulación de las matemáticas del primer ciclo universitario de ciencias: Estudio de la dinámica de poblaciones. Dans L. Ruiz-Higueras, A. Estepa & F. J. García (Éds), *Sociedad, escuela y matemáticas. Aportaciones de la teoría antropológica de lo didáctico (TAD)* (pp. 573-594). Jaén, Espagne : Publicaciones de la Universidad de Jaén.

Barquero, B., Bosch, M. & Gascón, J. (2011). Ecología de la modelización matemática: Restricciones transpositivas en las instituciones universitarias. Dans M. Bosch et al. (Éds), *Un panorama de la TAD* (pp. 553-577). Barcelone, Espagne : CRM.

Boas, R. (1981). Who needs those mean-value theorems, anyway ? *The Two-Year College Mathematics Journal*, *12*(3), 178-181.

Bonnet, P.-O. (1844). Démonstration simple du théorème de Fourier. *Nouvelles Annales de Mathématiques*, 1$^{re}$ série, tome 3, 119-121. http://archive.numdam.org/ARCHIVE/NAM/NAM_1844_1_3_/NAM_1844_1_3__119_1/NAM_1844_1_3__119_1.pdf

Bosch, M. & Gascón, J. (2001). *Las prácticas docentes del profesor de matemáticas.* http://www.ugr.es/~jgodino/siidm/almeria/Practicas_docentes.PDF

Bourdieu, P. & Passeron, J.-C. (1970). *La reproduction : éléments d'une théorie du système d'enseignement*. Paris : Les Éditions de Minuit.

Bourdieu, P. & de Saint-Martin, M. (1975). Les catégories de l'entendement professoral. *Actes de la recherche en sciences sociales*, *1*(3), 68-93.

Bourdieu, P. (1982). Les rites comme actes d'institution. *Actes de la recherche en sciences sociales*, *43*, 58-63.

Bourdieu, P. (2001). *Science de la science et réflexivité*. Paris : Raisons d'agir.

Chevallard, Y. (1991). *La transposition didactique. Du savoir savant au savoir enseigné* (2$^e$ éd.). Grenoble : La pensée sauvage.

Chevallard, Y. (2005). La place des mathématiques vivantes dans l'éducation secondaire : transposition didactique des mathématiques et nouvelle épistémologie scolaire. Dans C. Ducourtioux & P.-L. Hennequin (Éds), *La place des mathématiques vivantes dans l'éducation secondaire* (pp. 239-263). Paris : APMEP et Animath.

Chevallard, Y. (2011). *Journal du Séminaire TAD/IDD de l'année 2010-2011*. http://yves.chevallard.free.fr/spip/spip/IMG/pdf/journal-tad-idd-2010-2011-4.pdf

Daubelcour, J. P. (2009). *Évolution des programmes d'analyse et de géométrie au vingtième siècle en terminale scientifique*. http://jpdaubelcour.pagesperso-orange.fr/histoire20.html

Dieudonné, J. (1969). *Foundations of Modern Analysis*. New-York et Londres : Academic Press.

Gascón, J. (2001). Incidencia del modelo epistemológico de las matematicas sobre las prácticas docentes. *Revista Latinoamericana de Investigación en Matemática Educativa*, *4*(2), 129-159.

Gispert-Sambaz, H. (1982). *Camille Jordan et les fondements de l'analyse* (Thèse de 3$^e$ cycle).
http://www.maths.ed.ac.uk/~aar/jordan/gispert.pdf

Houzel, C. (1996). *Analyse mathématique, Cours et exercices, 1$^{er}$ cycle des Universités* (collection Sciences Sup). Paris : Belin.

Lombardi, H. (1999). À propos du théorème des accroissements finis. *Repères-IREM*, *34*, 55-69.

Ministère de l'Éducation nationale (1986). Programme de mathématiques de la classe de Terminale C. *Bulletin Officiel n° 31 du 11 septembre 1986.*

Ministère de l'Éducation nationale (2003). Programmes des CPGE mathématiques classe de première année MPSI. *Bulletin Officiel hors série n° 5 du 28 août 2003*, *Annexe 1*, 1154-55.

Monier, J.-M. (1999). *Analyse I, Cours et 300 exercices corrigés, 1$^{re}$ année MPSI, PCSI, PTSI*. Paris : Dunod.

Nelsen, R. B. (1993). *Proofs without words, exercises in visual thinking*. Washington : The Mathematical Association of America.

Ovaert, J.-L. & Verley, J.-L. (1983). *Analyse Vol. 1*, collection Léonard Épistemon. Paris : Cédic/Fernand Nathan.

Parra, V. E. & Otero, M. R. (2011). Praxeologías didácticas en la universidad y el fenómeno del «encierro»: Un estudio de caso relativo al límite y continuidad de funciones. Dans M. Bosch et al. (Éds), *Un panorama de la TAD. An overview of ATD* (pp. 719-741). Barcelone, Espagne : CRM.

Perrin, D. (1997). *Rigueur et formalisme(s)*. Dans M. Bailleul, C. Comiti, J.-L. Dorier, J.-B. Lagrange, B. Parzysz, M.-H. Salin (Éds), *Actes de la IX$^e$ école d'été de didactique des mathématiques* (Houlgate, 19-27 août 1997). Caen : ARDM et IUFM.

Rolle, M. (1690). *Traité d'Algèbre ou Príncipes généraux pour résoudre les questions de mathématique*. Paris : Estienne Michallet.

Shain, J. (1937). The method of cascades. *The American Mathematical Monthly*, *44*(1), 24-29.

Smith, D. E. (1984). *A source book in mathematics*. New York, NY : Dover. (Édition originale 1929)

Thurston, W. P. (1995). Preuve et progrès en mathématiques. *Repères IREM*, *21*, 7-26.

Tucker, T. W. (1997). Rethinking rigor in calculus: The role of the mean value theorem. *The American Mathematical Monthly*, *104*(3), 231-240.

Wittgenstein, L. (2005). *Recherches philosophiques.* Paris : Gallimard.